\documentclass{article}

\usepackage[english]{babel}
\usepackage{mathrsfs}

\usepackage{graphicx}
\usepackage{amsmath}
\usepackage{amssymb}
\usepackage{amsfonts}
\usepackage{color}
\usepackage{graphicx}

\usepackage{mathrsfs} %
\usepackage{graphicx}

\newtheorem{lemma}{Lemma}
\newtheorem{remark}{Remark}
\newtheorem{example}{Example}

\def\const{\mathrm{\rm const}}

\begin{document}
\large

\title{A Note on the Legendre Transformation}

\author{A.\,O. Remizov}
\date{}

\maketitle

\begin{abstract}
We present the Legendre transformation in a geometric way based on the procedure of the Legendrian lift.
This approach allows us to understand some interesting properties of it, in particular, the reason for the appearance of singularities of dual curves.
Also we consider application of the Legendre transformation to the Clairaut differential equation.
Finally, we say a few words class of contact transformations and present an infinite group of
contact transformations different from the Legendre transformation.

\end{abstract}

\section*{Introduction}

The paper has several goals. First, we present the Legendre transformation in a geometric way based on the procedure of the Legendrian lift.
Therefore, we start with the Legendre transformation of planar curves and then go to the Legendre transformation of functions. As we shall see below, this approach has some important advantages in compare with the standard definition of the Legendre transformation of functions, which is thrown in the reader's face without motivations.
Moreover, the presented geometric approach allows us to understand the meaning of restrictions in the standard definition of the Legendre transformation of functions
(for example, the convexity) and the reason for the appearance of singularities of dual curves.

Second, we consider application of the Legendre transformation to ordinary differential equations at the example of the Clairaut equation.
For this, the geometric interpretation of the Legendre transformation is more suitable.

Finally, we embed the Legendre transformation into a class of contact transformations and discuss some examples of contact transformations different from it.
In particular, we consider the pedal transformation (apparently the earliest example of contact transformations) and with its aid we obtain an infinite group
of contact transformations, which was first constructed by Sophus Lie. It should be remarked that many of results presented in the paper were previously known, but we presented these results from a more general point of view and explained the relationship between them.

We start with a small section presented some notions and results about planar curves and their singularities, which are necessary for understanding the main ideas of the paper.

\section{Curves and singularities}

Let $\gamma$ be an arbitrary curve on the plane given in the form
\begin{equation}
x = \phi(t),  \quad y = \psi(t),
\label{E1}
\end{equation}
where $\phi, \psi$ are smooth ($C^\infty$) functions.
Further, we shall assume that all considered functions and mappings are smooth ($C^\infty$) unless otherwise stated.

A point of the curve $\gamma$ corresponding the parameter $t_0$ is called  {\it singular} or {\it critical} if
\begin{equation}
\phi'(t_0) = \psi'(t_0) = 0.
\label{E-9Oct2021}
\end{equation}
Otherwise a point is called {\it regular}.
The curve $\gamma$ is called {\it regular} if all its points are regular.
In a neighborhood of every regular point, the curve is diffeomorphic to its tangent line, while for singular points it is not true.

\begin{figure}[!htp]
\centering
\includegraphics[scale=0.3]{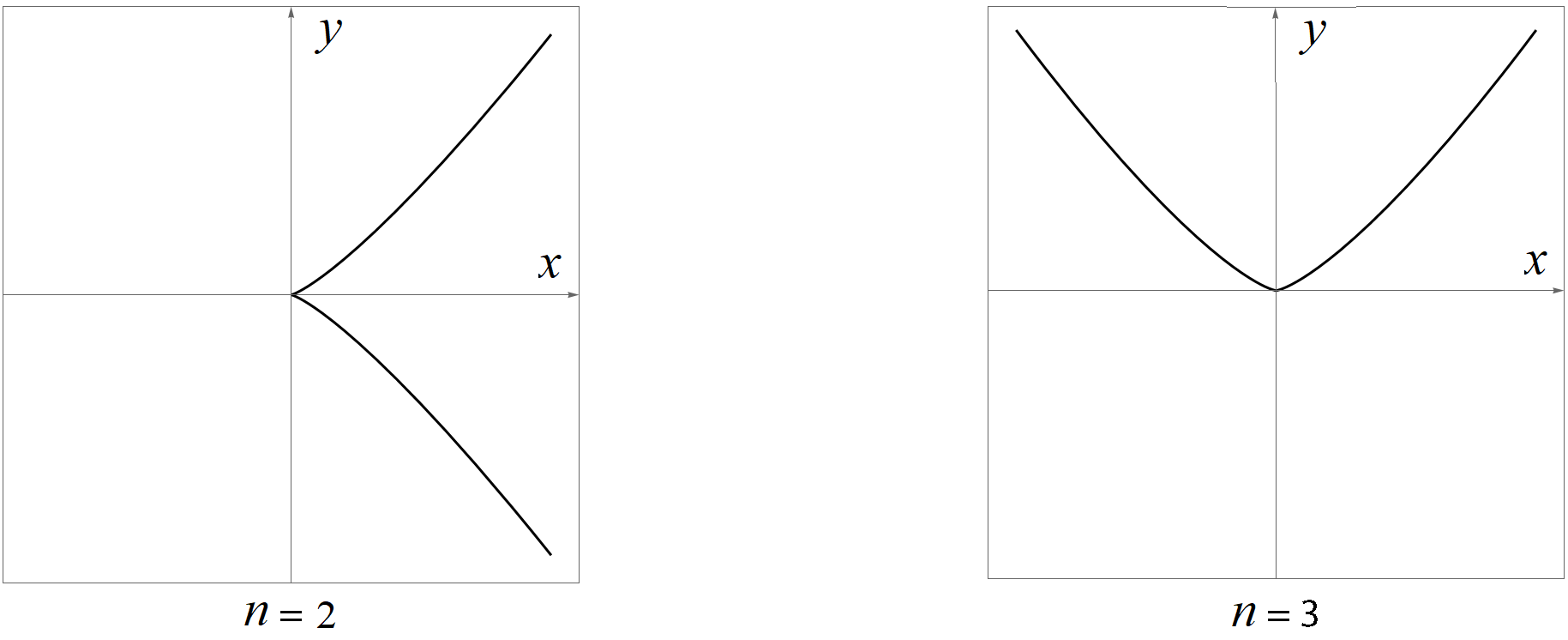}
\caption{\small
The curves \eqref{ide-29} with even $n$ (left) and odd $n$ (right).
}
\label{pic-ide7}
\end{figure}

\begin{example}
{\rm
In various applications, the following types of singular points often appear:
\begin{equation}
x = \phi(t) = \alpha t^n + o(t^n), \ \ y = \psi(t) = \beta t^{n+1} + o(t^{n+1}), \ \ \  \alpha \beta \neq 0.
\label{ide-29}
\end{equation}
Fig.~\ref{pic-ide7} presents such curves for even and for odd $n$.
The singular point of the curve \eqref{ide-29} with $n=2$ is called {\it semicubic cusp} or simply {\it cusp}.

It is proved that for $n=2$ or $3$ the germ of any curve \eqref{ide-29} at zero can be brought to the simplest form $x=t^n$, $y=t^{n+1}$
by means of an appropriate change of variables $(x,y)$ and a change of the parameter $t$.
However, for $n \ge 4$ it is not true~\cite{BG}.
}
\end{example}

\section{Legendre transformation}

The geometric definition of the Legendre transformation requires the notion of the Legendrian lift.
We start with the latter notion, which is an important mathematical procedure.
It can be applied to curves, functions and many other objects.
The Legendrian lift appears not only in mathematics and physics, but also in various technical applications and even in the living nature.
We shall say a few words about it below.

\subsection{Legendrian lift}

{\it The lift} or {\it $1$-graph} of the curve $\gamma$ given by formula \eqref{E1} is the curve
\begin{equation}
{\overline \gamma} : \quad x=\phi(t), \quad y=\psi(t),  \quad p=\frac{\psi'(t)}{\phi'(t)}
\label{E13}
\end{equation}
in 3-dimensional space with the coordinates $(x,y,p)$. This space denoted by $J^1$ is called the space of 1-jets of functions $y(x)$.

\begin{figure}[!htp]
\centering
\includegraphics[scale=1.0]{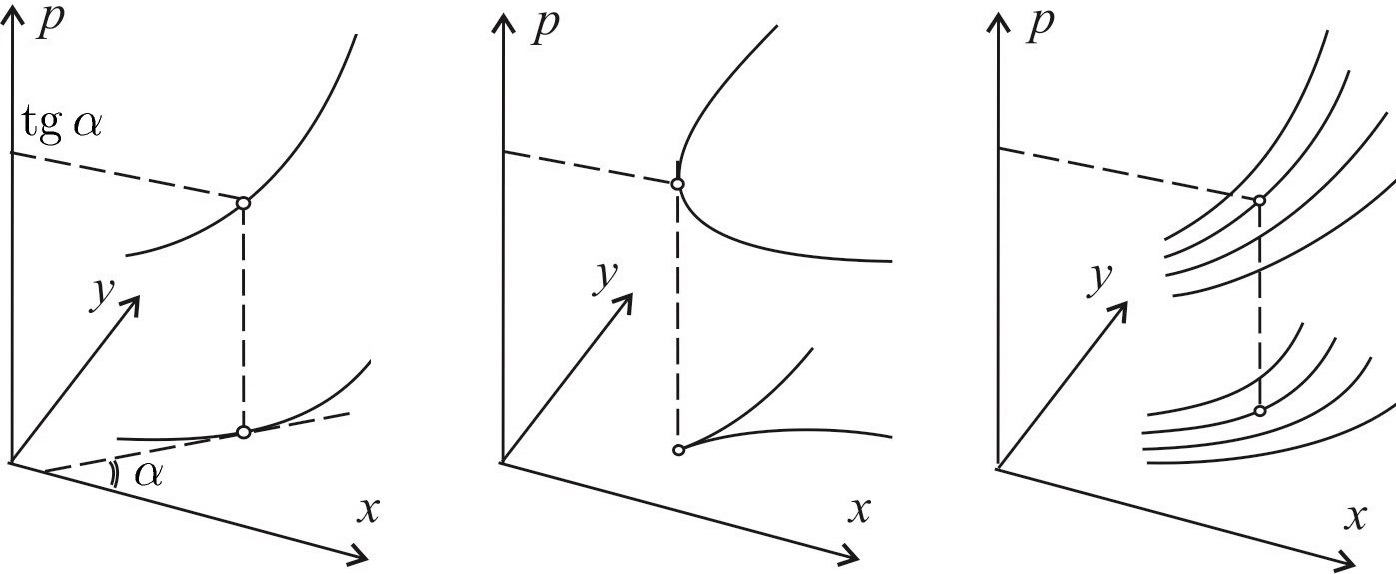}
\caption{\small
Legendrian lift of planar curves.
}
\label{ris06-4}
\end{figure}

In the latter equality in \eqref{E13} the denominator $\phi'(t)$ may vanish, and we have to discuss the correctness of such operation.
If $\phi'(t) = 0$, but $\psi'(t) \neq 0$, we have $p=\infty$ at the considered point, that is, the tangent direction of the curve $\gamma$  is parallel to the $y$-axis.
The dividing by zero disappears, if we interchange the $x$-axis with the $y$-axis.
The situation is more complicated if the curve $\gamma$ has a singular point where $\phi'(t) = \psi'(t) = 0$.
However, in the most of interesting cases the function ${\psi'(t)}/{\phi'(t)}$ has the same one-sided limits,
and after an appropriate definition at the point itself, the function ${\psi'(t)}/{\phi'(t)}$ becomes continuous.
For example, such a situation takes place for the semicubic parabola, see Fig.~\ref{ris06-4}.

\begin{remark}
{\rm
The $1$-graph of the curve $\gamma$ does not depend on the choice of the parametrization the curve $\gamma$.
}
\end{remark}

According to modern biological concepts going back to the famous work of Hubel and Wiesel (1959), the brain of mammals performs the Legendrian lift of apparent contours of flat images perceived by the retina. There exists an area of the primary cortex of the brain associated with vision (V1) that contains neurons united in groups with complex hierarchy (columns, hypercolumn). Each column is sensitive to a certain point of the retina and a certain tangent direction at this point.
Flat images perceived by the retina are lifted in the space of 1-jets, and further work of the brain occurs with the lifted objects.
See \cite{Peti}  and the references therein.
The only difference between the lift performed by the brain and the lift described above is that the brain uses the third coordinate $\theta = \arctan p$ instead of $p$.
It is very natural: $\theta$ is more appropriate for procedure performed biologically, since it is bounded, while $p$ is not.

The theoretical concept presented above has practical applications: it is the basis of some inpainting algorithms that allow to recover corrupted images.
In these algorithms, the Legendrian lift also plays a crucial role:
First, we lift the image, which is considered as a family of contours, in 3-dimensional space and then we act to the lifted object a hypoelliptic diffusion
associated with the contact structure. For details, see  \cite{Bos1, Bos2}.

\subsection{Legendre transformation and duality}

{\it Legendre transformation} is an automorphism $\Lambda: J^1 \to J^1$ given by the formula $(x,y,p) \mapsto (X,Y,P)$,
where the coordinates are connected by the equalities:
\begin{equation}
x = P, \quad  p = X, \quad y+Y = xp = XP.
\label{E14}
\end{equation}

The Legendre transformation can be applied to various objects: curves, functions, differential equations, etc.
The objects obtained in this way are called {\it dual} of the original, they are often denoted by the same symbols with asterisk.
Dual objects reflect many important properties of their preimages.

\begin{figure}[!htp]
\centering
\includegraphics[scale=0.5]{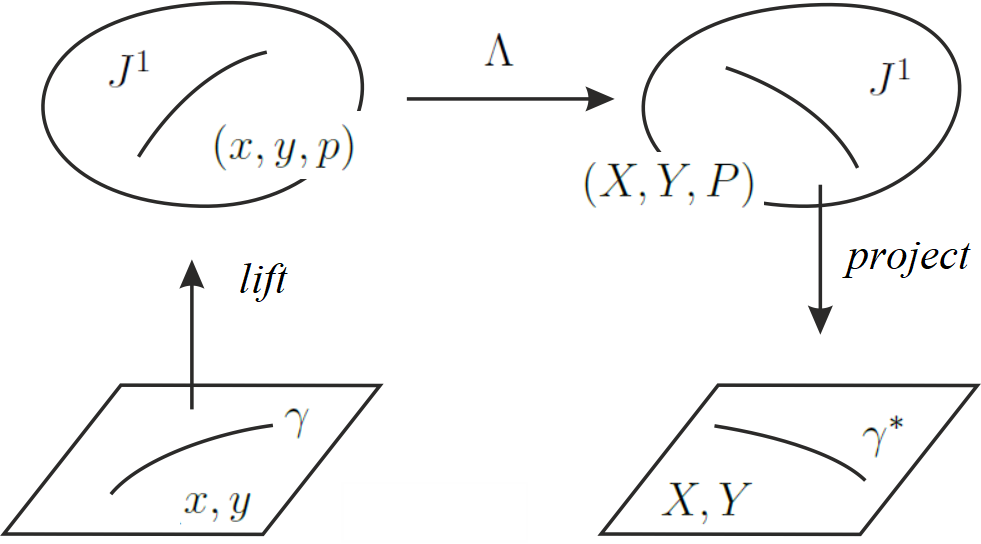}
\caption{\small
Definition of the Legendre transformation of a planar curve.
}
\label{Legen}
\end{figure}

Consider the application of the Legendre transformation to the curve $\gamma$  given by formula \eqref{E1}.
It is expressed by the diagram presented in Fig.~\ref{Legen}.

In Fig.~\ref{Legen}, the vertical upward arrow denotes the Legendrian lift of the curve $\gamma$ to the space $J^1$ with coordinates $(x,y,p)$.
The horizontal arrow denotes the Legendre transformation $\Lambda: J^1 \to J^1$ given by formula \eqref{E14},
the vertical downward arrow denotes the projection $\pi: (X,Y,P) \mapsto (X,Y)$.
The composition of these three mappings results in the Legendre transformation of the curve $\gamma$ on the $(x,y)$-plane to a certain curve $\gamma^*$ on the $(X,Y)$-plane,
which is called {\it dual} to $\gamma$.

The curve $\gamma^*$ dual to a curve $\gamma$ can have singular points even $\gamma$ is regular.
Obviously, they appear due to the projection $\pi: (X,Y,P) \mapsto (X,Y)$.

\begin{example}
\label{Ex-Dual-1}
{\rm
Let us find the dual curves for the line $\gamma_1: y=at+b$, the parabola $\gamma_2: y=x^2$, and the cubic parabola $\gamma_3: y=x^3$.

The line $\gamma_1$ can be given in the parametric form $x=t$, $y=at+b$.
Hence $p=a$ and, consequently, we get $X=a$, $Y=-b$, that is, the dual curve $\gamma_1^*$ is a single point on the plane.

The parabola $\gamma_2$ can be presented in the parametric form $x=t$, $y=t^2$.
Hence $p=2t$ and, consequently, we get $X=2t$ and $Y=t^2$.
Thus, the dual curve $\gamma_2^*$ is the parabola $Y=X^2/4$.

The cubic parabola $\gamma_3$ can be presented in the parametric form $x=t$, $y=t^3$.
This yields $p=3t^2$ and $X=3t^2$, $Y=2t^3$.
This shows that the dual curve $\gamma_3^*$ is the semicubic parabola with the cusp at the origin.
The cusp corresponds to the inflection point of the initial curve~$\gamma_3$.
}
\end{example}

\begin{lemma}
\label{Lem-1}
The Legendre transformation has the following properties.

1. The curve $\gamma^*$ dual to the curve $\gamma$ of the form \eqref{E1} is regular at all points that correspond to points of the curve $\gamma$
with non-zero curvature $\kappa(t) \neq 0$. The curve $\gamma^*$ has simicubic cusp at points that correspond to points of $\gamma$ where
$\kappa(t) = 0$, but $\kappa'(t) \neq 0$.

2. Assume that the curve $\gamma$ is the graph of a functions $y=f(x)$ such that
\begin{equation}
f(0) = \frac{df}{dx}(0) = \cdots = \frac{d^{n}f}{dx^{n}}(0)=0, \ \ \, \frac{d^{n+1}f}{dx^{n+1}}(0) \neq 0.
\label{E5-0}
\end{equation}
Then, in the neighborhood of the origin, the curve $\gamma^*$ dual to $\gamma$ has the form \eqref{ide-29}.

3. The Legendre transformation $\Lambda$ is an {\it involution}, that is, $\Lambda^2$ is the identity transformation.
This is equivalent to the relation $\Lambda = \Lambda^{-1}$, which yields $(\gamma^*)^* = \gamma$.
\end{lemma}

The proof is by direct calculation.

\medskip

The lemma shows that if $\gamma$ is the graph of a function $y=f(x)$, then
regular points of $\gamma^*$ correspond to points where $f''(x) \neq 0$, while
semicubic cusps correspond to points where
$$
f''(x) = 0, \ \ f'''(x) \neq 0,
$$
i.e., simple (cubic) inflection points of the function $f$.

\begin{example}
{\rm
In Fig.~\ref{ris06-5}, the curve dual to the graph of sine is presented.
(For better visibility, the scale along the $x$-axis is increased approximately 20~times.)
Cusps of the dual curve correspond to points $x = \pi n$ with integer $n$.

\begin{figure}[!htp]
\centering
\includegraphics[scale=0.25]{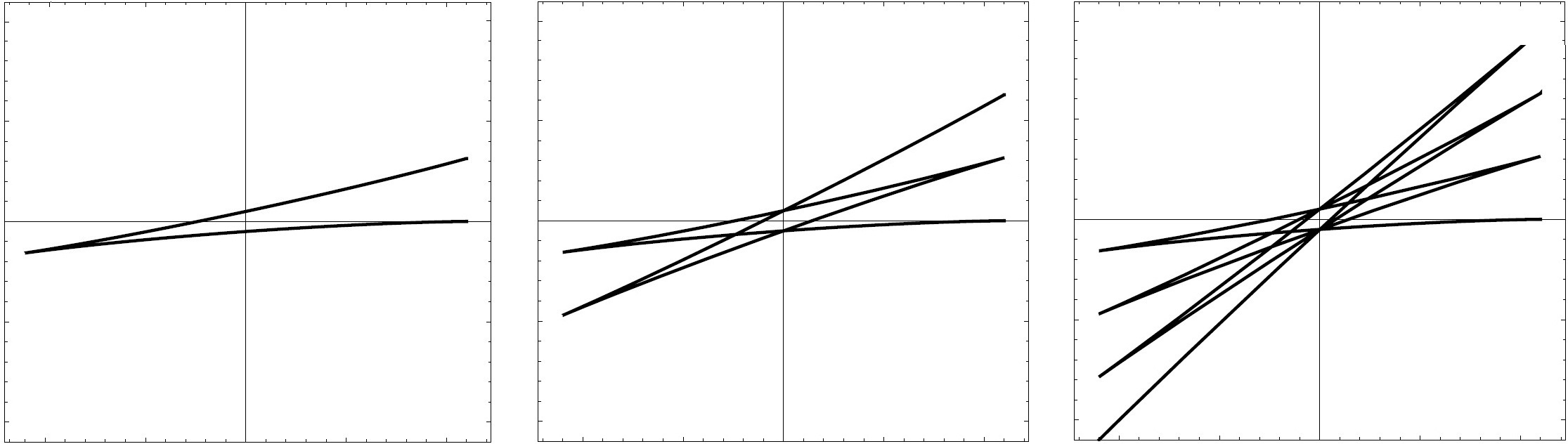}
\caption{\small
The curve dual to $y=\sin x$.
From the left to right: $x$ runs from $0$ to $2\pi$, $4\pi$, $7\pi$, respectively.
}
\label{ris06-5}
\end{figure}
}
\end{example}

\begin{remark}
{\rm
It is worth observing that there exist other ways to define dual curves.
For instance, in algebraic geometry the definition of dual curves is based on the same construction, where the mapping $\Lambda$ given by formula~\eqref{E14}
is replaced with the mapping $(x,y,p) \mapsto ({\bar X}, {\bar Y}, {\bar P})$, where
\begin{equation}
{\bar X} = X/Y, \ \ {\bar Y} = 1/Y, \ \ {\bar P} = d{\bar Y}/d{\bar X},
\label{E15}
\end{equation}
$X,Y$ are defined in \eqref{E14}.
}
\end{remark}

\begin{example}
{\rm
Consider the curve dual to the ellipse $(x/a)^2+(y/b)^2=1$. It is easy to see that this curve is the hyperbola $Y^2-(aX)^2=b^2$ in the coordinates $(X,Y)$ and
the ellipse $(a {\bar X})^2 + (b {\bar Y})^2 = 1$ in the coordinates $({\bar X},{\bar Y})$.
It is not surprising, because all conic sections (which include ellipses and hyperbolas) are projectively equivalent and
the coordinates $(X,Y)$ and $({\bar X}, {\bar Y})$ are connected by a projective transformation \eqref{E15}.
}
\end{example}

Now consider an important case that the curve $\gamma$ is the graph of a function $y=f(x)$.
It is naturally to ask: whether the dual curve $\gamma^*$ is the graph of a certain function $Y=f^*(X)$ as well.
If so, the function $f^*$ is called the dual function for $f$ or the Legendre transformation of $f$.

Examples considered above show that the dual curve for the graph of a function is not always the graph of a function.
Lemma~\ref{Lem-1} gives the following sufficient condition:
if the function $f$ is strictly convex ($f''>0$) or $f$ is strictly concave ($f''<0$) on its domain, then
the dual curve is the graph of a function, whence the dual function $f^*$ exists.
However, the function $f^*(X)$ can be defined not for all $X$.
Indeed, since $X=p$, the domain of $f^*$ coincides with the range of $f'(x)$, and $f^*$ is defined for $X$ that belong to the range of $f'$ only.

\begin{figure}[!htp]
\centering
\includegraphics[scale=1.1]{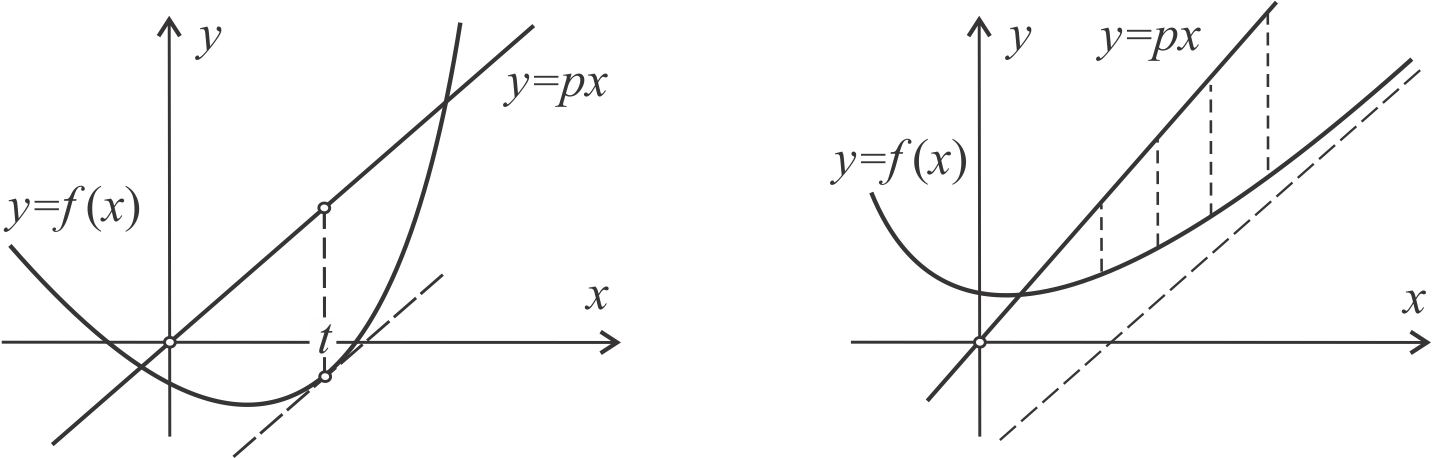}
\caption{\small
Illustration of formula \eqref{E-sup}. On the right: the case that $p \notin R$.
}
\label{function}
\end{figure}

If the function $f$ is strictly convex, its Legendre transformation $f^*$ can be defined by the formula
\begin{equation}
f^*(p) = \sup\limits_{x} (xp-f(x)), \ \  p \in R,
\label{E-sup}
\end{equation}
see Fig.~\ref{function}.
If the function $f$ is strictly concave, its Legendre transformation $f^*$ can be defined by the formula
\begin{equation}
f^*(p) = \inf\limits_{x} (xp-f(x)), \ \  p \in R,
\label{E-sup2}
\end{equation}
where in the both formulas $R$ denotes the range of $f'$.

Let us prove formula \eqref{E-sup} (for \eqref{E-sup2} the proof is similar).
Fixing arbitrary $p \in R$, consider the function $F(x):= xp - f(x)$.
Since the function $f$ is strictly convex, $F$ is strictly convex as well. Therefore, $F$ has no more than one critical point.

Assume that $F$ has a critical point, denote it by $t$ (Fig.~\ref{function} left).
Then $F'(t)=0$ and $f'(t)=p$, whence from \eqref{E-sup} we get:
\begin{equation}
\label{E-sup3}
f^*(p) = \sup\limits_{x}\, (xp-f(x)) = tp -f(t) = tf'(t) - f(t), \ \ p = f'(t).
\end{equation}
The graph of the function $f^*(p)$ given by formula \eqref{E-sup3} coincides with the Legendre transformation of the graph $y=f(x)$.
Indeed, let us introduce the parameter $t=x$, then $y=f(t)$, $p=f'(t)$, This yields
$$
X = f'(t), \ \ \, Y = tf'(t) - f(t),
$$
that is, \eqref{E-sup3} is equivalent to $f^*(X) = Y$.

If for the chosen $p$ the function $F$ has no critical point (Fig.~\ref{function} right),
the value $f^*(p)$ cannot be defined by either the first or second way.
However, in such case one can set $f^*(p) = +\infty$ according to~\eqref{E-sup}.

\begin{remark}
{\rm
Formulas \eqref{E-sup} and \eqref{E-sup2} are often used as the definition of the Legendre transformation of functions.
In our opinion, such a definition has a serious drawback compared to the geometric approach presented above: it is not natural and well motivated.
}
\end{remark}

\begin{remark}
{\rm
The Legendre transformation can be obviously adapted for hypersurfaces and, consequently, for functions depending on any number of variables.
In particular, this way brings us to the same formulas \eqref{E-sup} and \eqref{E-sup2}, where
$x$ and $p$ belong to the spaces of the same dimension and $xp$ means the dot product in orthonormal coordinates.
}
\end{remark}

\subsection{Clairaut differential equation}

The Legendre transformation can be applied to differential equations as well.

Consider the first-order differential equation
\begin{equation}
F(x,y,p)=0, \ \ p=dy/dx.
\label{E-diffur}
\end{equation}
Passing to the new variables $X,Y,P$ by formula \eqref{E14} and substituting the corresponding expressions in \eqref{E-diffur},
we obtain a new first-order differential equation, which is called {\it dual} to the initial equation \eqref{E-diffur}:
\begin{equation}
\label{E-diffur2}
F^*(X,Y,P) = 0, \ \ P=dY/dX,
\end{equation}
where
$$
F^*(X,Y,P) = F(P,XP-Y,X).
$$
In fact, here we have to verify the equality $P=dY/dX$ only, which follows from the equality $p=dy/dx$ and formula \eqref{E14}.
It more convenient to write the both equalities via differential 1-forms and check that 1-forms $PdX-dY$ and $pdx-dy$ vanish simultaneously:
$$
PdX-dY = xdp - d(xp-y) = xdp - d(xp) + dy = - (pdx - dy).
$$

From geometric viewpoint, this means that the Legendre transformation \eqref{E14} sends {\it contact planes} $pdx-dy=0$ to {contact planes} $PdX-dY=0$, i.e.,
it preserves the natural contact structure of the space~$J^1$. This implies that the Legendre transformation sends integral curves of equation \eqref{E-diffur}
to integral curves of its dual equation \eqref{E-diffur2}. In other words, integral curves of equation \eqref{E-diffur2} are dual to integral curves of equation \eqref{E-diffur}, and vice versa (since the Legendre transformation is involution).
This fact can be used for solution of differential equation, if the dual equation can be solved easier that the initial one.
The notion of duality also helps to study singularities of implicit differential equations; see \cite{BT}.

A nice example of application the Legendre transformation to differential equation is the {\it Clairaut equation}, which has the form
\begin{equation}
xp - y = f(p), \ \ p=dy/dx,
\label{E16}
\end{equation}
where $f$ is an arbitrary function. See, for example, \cite{Arnold, MK}.

The Legendre transformation sends differential equation \eqref{E16} to the equation $Y=f(X)$, which does not contain the derivative.
Therefore, the Legendre transformation sends integral curves of equation \eqref{E16} to points of the $(X,Y)$-plane filling the graph of the function $Y=f(X)$.
Taking into account Example~\ref{Ex-Dual-1}, one can see that integral curves of equation \eqref{E16} are straight lines tangent to a curve dual to the graph $Y=f(X)$.

A direct calculation shows that this curve given by formula
\begin{equation}
x = f'(p), \ \ y = pf'(p)-f(p)
\label{E17}
\end{equation}
is the {\it discriminant curve} of equation \eqref{E16} and the {\it envelope} of the family of solutions of equation \eqref{E16} -- its tangent lines.
This gives a simple way to get all solutions of equation \eqref{E16}:
First, we find the discriminant curve \eqref{E17}, which is called {\it singular solution} of \eqref{E16}.
The remaining solutions are tangent lines to the discriminant curve taken at all its points.
Moreover, it is not hard to see that these tangent lines  are given by the formula
\begin{equation}
y= cx-f(c), \ \ c = {\rm const}.
\label{E17-}
\end{equation}
Remark that formula \eqref{E17-} can be obtained from equation \eqref{E16} after replacing $p$ with $c$.

\begin{figure}[!htp]
\centering
\includegraphics[scale=0.65]{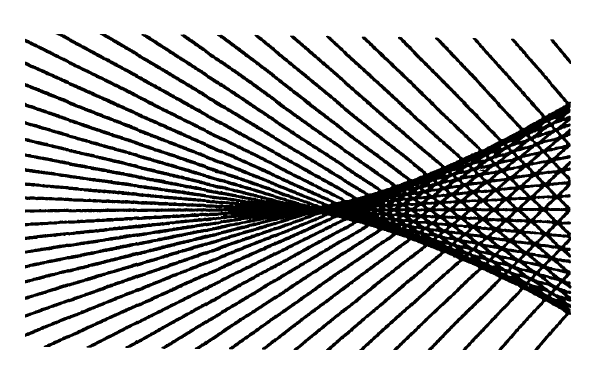}
\caption{\small
Integral curves of the Clairaut equation $xp - y = p^3$ are tangent lines to the semicubic parabola~\eqref{semic}.
The picture is taken from~\cite{IzumiyaKurokawa1}
}
\label{ris06-6}
\end{figure}

\begin{example}
{\rm
In Fig.~\ref{ris06-6}, integral curves of equation \eqref{E16} with the function $f(p)=p^3$ are presented.
In this case, the discriminant curve \eqref{E17} is the semicubic parabola
\begin{equation}
x=3p^2, \ \  y=2p^3
\label{semic}
\end{equation}
with the cusp at the origin.
The family of tangent lines to the curve \eqref{semic} has concentration (patch) near this curve and especially near its cusp.

This phenomenon is called a {\it caustic} \cite{Arnold2}.
In physics, a caustic is the envelope of light rays which have been reflected or refracted by a curve (on the plane) or a surface (in the space). Fig.~\ref{ris06-6} contains a caustic~-- the semicubic parabola~\eqref{semic}.
}
\end{example}

Finally, let us consider the Clairaut equation of more general type:
\begin{equation}
F(p, xp - y) = 0, \ \ p=dy/dx,
\label{E16001}
\end{equation}
where $F(u,v)$ is an arbitrary function of two variables.
The Legendre transformation sends \eqref{E16001} to the equation $F(X,Y)=0$ without derivative.
Proceeding in the same way as before, we conclude that integral curves of equation \eqref{E16001} are straight lines tangent to
the curve dual to $F(X,Y)=0$. The equation $F(X,Y)=0$ defines the discriminant curves of equation \eqref{E16001}, which is its singular solution
and the envelope of other solutions~-- tangent lines given by the formula
$$
y=ax+b, \ \  F(a,-b)=0,
$$
which is a generalization of~\eqref{E17-}.

\begin{example}
{\rm
Consider the equation
\begin{equation}
p^3 = (xp-y)^2, \ \ p=dy/dx,
\label{E16002}
\end{equation}
which has the form \eqref{E16001} with $F(u,v) = u^3-v^2$.

The curve dual to the curve $u^3-v^2=0$ is the cubic parabola $y = \tfrac{4}{27}x^3$, which is the discriminant curve of equation~\eqref{E16002}
and its singular solutions. The family of non-singular solutions of \eqref{E16002} consists of all tangent lines to the singular solution, which are given by the formula
$$
y=ax+b, \ \  a^3 = b^2.
$$
The latter equation can be written in the form
\begin{equation}
a^3 = (ax-y)^2, \ \ a = \const,
\label{Wright}
\end{equation}
or in the form $y = c^2 x + c^3$, where $a = c^2$.
Formula \eqref{Wright} can be obtained from equation \eqref{E16002} by replacing the derivative $p$ with the constant $c$.
}
\end{example}

\section{Contact transformations}

The Legendre transformation is not a unique transformation $J^1 \to J^1$ preserving the contact structure of this space.
There exist a large group of diffeomorphisms of the space that also preserve the contact structure, they are called {\it contact transformations}.
Equivalent definition: contact transformations are transformations of curves in the plane in which tangent curves are transformed into tangent curves.

Let us try to determine condition under which a diffeomorphism $(x,y,p) \to (X,Y,P)$ (global or local) preserves the contact structure of the space~$J^1$.
Namely, consider a diffeomorphism $(x,y,p) \to (X,Y,P)$ given by the formula
\begin{equation}
X = F(x,y,p), \ \ Y = G(x,y,p), \ \ P = H(x,y,p),
\label{APP5-000}
\end{equation}
where the functions $F, G$ satisfy some relation (that we shall find) and the function $H$ is uniquely defined by $F, G$ due to the condition $P=dY/dX$.
Therefore, the condition of preservation the contact structure has the form
\begin{equation}
H(x,y,p) = \frac{dY}{dX} = \frac{G_xdx+G_ydy+G_pdp}{F_xdx+F_ydy+F_pdp} = \frac{G_x+pG_y+p'G_p}{F_x+pF_y+p'F_p},
\label{APP5-6}
\end{equation}
where $p' = {dp}/{dx}$.

Since $P = H(x,y,p)$ does not depend on $p'$, the right hand side of equality \eqref{APP5-6} does not depend on $p'$ as well.
This yields
$$
\frac{G_x+pG_y}{F_x+pF_y} \equiv \frac{G_p}{F_p}
$$
or, equivalently,
\begin{equation}
{F_p} (G_x+pG_y) \equiv {G_p}(F_x+pF_y).
\label{APP5-7}
\end{equation}
The identity \eqref{APP5-7} connecting the functions $F$ and $G$ is a necessary and sufficient condition for the diffeomorphism \eqref{APP5-000} be a contact transformation.
It is not hard to see that the Legendre transformation satisfies \eqref{APP5-7}.

\medskip

The Legendre transform is the most commonly used and the most useful contact transformation, but historically it was not the first one.
Apparently, the first contact transformation was the so-called {\it pedal transformation}, which was investigated by Colin Maclaurin, Arthur Cayley, and Sophus Lie.

To define pedal transformation, let us fixe an arbitrary point $O$ on the plane (without loss of generality, $O$ can be the origin).
The point $O$ is called the {pole} or {center} of the transformation.
Let $\gamma$ be a curve on the plane. Then the {\it pedal curve} of $\gamma$ is the locus of points $X$ so that the line $OX$ is perpendicular to the tangent line to $\gamma$ passing through the point~$X$.

\begin{example}
{\rm
The pedal curve of the circle is the limacon (Pascal's snail), the special shape of which is determined by the distance $\rho$ between the pole $O$ and the circle (Fig.~\ref{circle-pedal}). Considering $\rho$ as a parameter, we get a family of the pedal curves, some of which have singular points.
Namely, when the pole $O$ lies on the circle itself ($\rho=0$), the pedal curve has the cusp at~$O$, this curve is called {\it cardioid} (the third figure from the left).

\begin{figure}[!htp]
\centering
\includegraphics[scale=1.1]{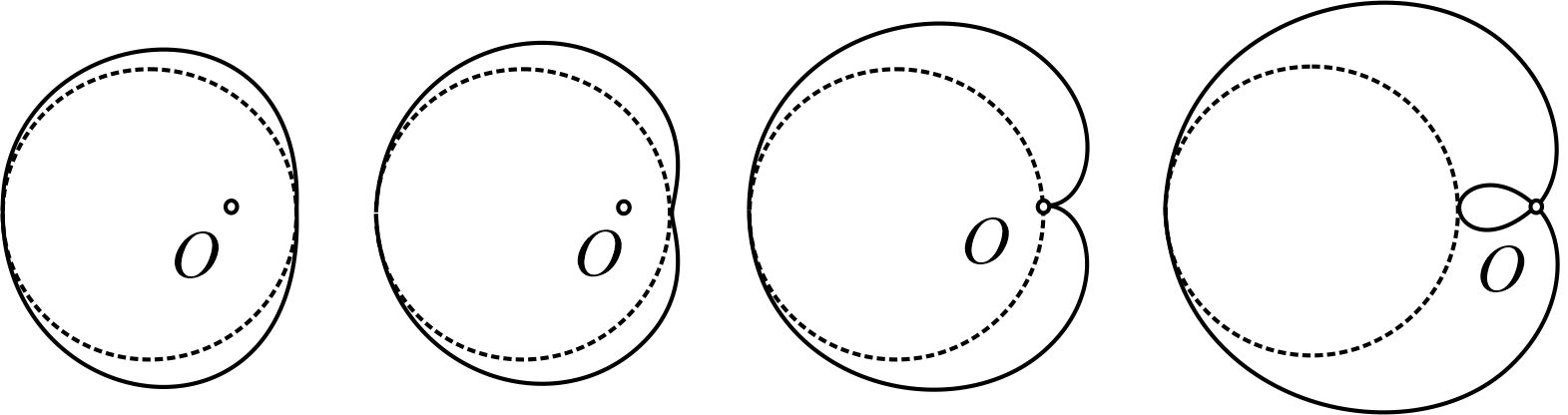}
\caption{\small
The pedal curves of the circle (depicted as the dotted line) are limacons, whose shape depends on
the distance between the pole $O$ and the circle
}
\label{circle-pedal}
\end{figure}

}
\end{example}

\begin{lemma}
\label{Lem-2}
In the Cartesian coordinate system centered at $O$, the pedal transformation is given by the formula
$$
{\mathcal P}: (x,y,p) \mapsto ({\bar X}, {\bar Y}, {\bar P}),
$$
where
\begin{equation}
{\bar X} = XY/(1+X^2), \ \ {\bar Y} = -Y/(1+X^2), \ \ {\bar P} = d{\bar Y}/d{\bar X},
\label{E15000}
\end{equation}
the relation between $X,Y, P$ and $x,y,p$ is given in \eqref{E14}.
\end{lemma}

The proof is by direct calculation.

\medskip

In contrast to the Legendre transformation $\Lambda$, the pedal transformation ${\mathcal P}$ is not an involution.
Therefore, one can consider the degrees (compositions) of the pedal transformation.
Fixing the pole $O$, consider the cyclic group generated by ${\mathcal P}$, that is, consisting of elements ${\mathcal P}^n$ with integer~$n$:
$$
\ldots \ {\mathcal P}^{-3}, \ {\mathcal P}^{-2}, \ {\mathcal P}^{-1}, \ {\mathcal P}^0, \ {\mathcal P}^1, \ {\mathcal P}^2, \ {\mathcal P}^{3} \ \ldots
$$
A natural question: is thus group finite or infinite? We will get an answer below.

\medskip

Sophus Lie suggested a nice way how to extend this discrete grout to a continuous group.
For this, consider the polar coordinates $(r, \phi)$ on the $(x,y)$-plane and the polar coordinates $(R, \Phi)$ on the $({\bar X},{\bar Y})$-plane.
As the third coordiante in the space $J^1$ let us take the difference between the angle of inclination of the radius vector to a given point on the plane
and the angle of inclination of the tangent to the curve at this point such that the orders of subtraction in the preimage and image are opposite:
\begin{equation}
\begin{aligned}
x  &=r\cos \phi, \ \ \ y=r\sin \phi, \ \ \ \ \psi = \phi - \arctan p,  \\
{\bar X} &=R\cos \Phi, \ \ {\bar Y} = R\sin \Phi, \ \ \Psi = \arctan {\bar P} - \Phi,
\end{aligned}
\label{APP5-4}
\end{equation}
where $x,y,p$ and ${\bar X}, {\bar Y}, {\bar P}$ are connected via~\eqref{E15000}.

It is not hard to check that in the coordinates \eqref{APP5-4}, the transformation ${\mathcal P}^n$ with every integer $n$ is given by the formula
\begin{equation}
R = r |\sin \psi|^n, \quad
\Phi = \phi - n(\psi + \tfrac{\pi}{2}), \quad
\Psi = \psi.
\label{APP5-5}
\end{equation}
The obtained formulas show shows that ${\mathcal P}^{n} \neq {\mathcal P}^{m}$ for any $n \neq m$.
Therefore, the cyclic group $\{{\mathcal P}^{n}\}$ is infinite.

Then let us consider formula \eqref{APP5-5} with real $n$ and the operation
${\mathcal P}^{n} \cdot {\mathcal P}^{m} = {\mathcal P}^{n+m}$.
This yields a continuous group of contact transformations (it was first constructed by Sophus Lie).

\medskip

We will not go deep into the theory of contact transformations,
referring the interested reader to the specialized literature, which is quite extensive.
For example, see \cite{Feld, Giz, Lie}.


\end{document}